# Towards Noncommutative Topological Quantum Field Theory: Tangential Hodge-Witten cohomology

I P ZOIS[1]
PPC, TRSC, 9, Leontariou Street,
Kantza, Pallini, GR-153 51, Athens, Greece
and
The American College of Greece, Agia Paraskevi, GR-153 42, Athens Greece
E-mail: i.zois@exeter.oxon.org

Abstract. A few years ago we initiated a program to define Noncommutative Topological Quantum Field Theory (NCTQFT for short, see [1]). The motivation came both from physics and mathematics: On the one hand, as far as physics is concerned, following the well-known holography principle of 't Hooft (which in turn appears essentially as a generalisation of the Hawking formula for black hole entropy), quantum gravity should be a topological quantum field theory. On the other hand as far as mathematics is concerned, the motivation came from the idea to replace the moduli space of flat connections with the Gabai moduli space of codim-1 taut foliations for 3 dim manifolds. In most cases the later is finite and much better behaved and one might use it to define some version of Donaldson-Floer homology which, hopefully, would be easier to compute. The use of foliations brings noncommutative geometry techniques immediately into the game. The basic tools are two: Cyclic cohomology of the corresponding foliation C*-algebra and the so called "tangential cohomology" of the foliation. A necessary step towards NCTQFT is to develop some sort of Hodge theory both for cyclic (and Hochschild) cohomology and for tangential cohomology. Here we present a method to develop a Hodge theory for tangential cohomology of foliations by mimicing Witten's approach to ordinary Morse theory by perturbations of the Laplacian.

## I. Introduction

Let F be a smooth p-dim foliation on a closed n-dim manifold M (hence $q = n - p$ is the codim), equipped with an invariant transverse measure $\Lambda$. It is well known that there exist real valued smooth functions on M having only Morse or birth - death singularities. We shall denote by h (resp v) the horizontal or tangential (resp vertical or transverse) local coordinates and by $L_x$ the leaf through the point x on M.

For any smooth real function on M we denote by $d_F$ the differential of in the leaf (horizontal or tangential) directions. A point a on M for which the leaf differential vanishes will be called a tangential singularity for . For such a singularity the horizontal or tangential Hessian $d^2_F$ makes sense and in local coordinates (h,v) one has

$$d_F = \sum_{1 \leq i \leq p} \partial/\partial h^i (h,v)$$

and

$$d^2_F (h,v) = ((\partial^2/\partial h^i \partial h^j)(h,v))_{ij}$$

---
[1] To whom any correspondence should be addressed.

The index of a tangential singularity a on M is defined as the number of minus signs in the signature of the quadratic form $d^2_F$ (a).

Definition *1:* A tangential singularity a on M of a smooth real function on M as above is called a Morse singularity if $d^2_F$ (a) is non-singular.

We denote by T( ) (resp M( ), $M_i$( )) the set of all tangential singularities (resp. of Morse singularities, Morse singularities of index i, where $0 \leq i \leq p$) of the function . The first complication emerges since in this case a good definition for a tangential (or horizontal) Morse function cannot be reduced to simply a smooth function on M having only tangential Morse singularities. This is explained in the following Lemma:

Lemma *1*: Let (M,F) be as above. Assume that there exists a smooth function again as above with only tangential Morse singularities. Then the set of all tangential Morse singularities of is a closed q-dim submanifold transverse to the foliation.

Proof: We suppose that is a smooth function on M such that for any leaf L in the quotient space M/F the restriction of on the leaf L has no degenerate critical points. Then the map

$$x \text{ maps to } (x, d_F(x))$$

from M to T*F is transverse to the zero section of T*F since its differential is given on any foliation chart $\varphi = U \times T$ by

$$d(d_F)(h,v)(X_h, X_v) = ((X_h, X_v), \partial_{hh}(h,v)X_h + \partial_{hv}(h,v)X_v)$$

where the subscript h, v denotes partial derivative with respect to the corresponding coordinates and $\det(\partial_{hh}(h,v)) \neq 0$ for a Morse singularity with coordinates (h,v). This implies first that the set of all Morse singularities M( ) is a closed submanifold of M with dim(M( )) = codim(F) and second that M( ) is transverse to the foliation F because for any non-zero tangent vector $X = (X_h, X_v)$ of M( ) at the point (h,v), one has that

$$\partial_{hh}(h,v)X_h + \partial_{hv}(h,v)X_v = 0.$$

This means that the transverse component $X_v$ of $X \neq 0$ is non-zero which proves that M( ) is transverse to the foliation F. This concludes the proof.

It turns out that many interesting foliations have no closed transversals and hence any good notion of tangential Morse function should allow degenerate critical points in the leaf direction. (However taut foliations which are the ones appearing in the Gabai moduli space [2] do have closed transversals).

*Definition 2:* We call almost Morse function a smooth function as above with degenerate critical points which only occure at a negligible set of leaves (namely we allow degenerate critical points but not too many).

*Definition 3:* A good almost Morse function is an almost Morse function which is generically unfolded in the sense of Igusa [3] (roughly this means that it has only birth-death singularities, namely points where critical points cancel or created in pairs).

## II. Witten's perturbation by a Morse function [4]-Tangential version

Let (M,F) be a foliation as above equipped with a holonomy invariant transverse measure $\Lambda$. We choose a smooth Riemannian metric on M and denote by $\Delta^k_L$ ($0 \leq k \leq p$) the corresponding Laplace operator on the leaf L acting on k-forms. We know that the bundle of Hilbert spaces is square integrable and thus has a well-defined Murray-von Neumann dimension

$$\beta_k = \dim_\Lambda (\text{Ker}(\Delta^k_L)) < \infty$$

which does not depend on the choice of metric. Assume moreover that codim(F) $\leq$ dim(F).

Let $\varphi$ be a smooth real function on M which is good almost Morse function and $\tau$ a positive real parameter. For each leaf L and $0 \leq k \leq p$ we denote by $d^k_{\tau,L}$ the closure (in $L^2(L, \wedge^* F)$, space of square integrable forms on the leaf L) of the operator which sends each smooth k-form $\alpha$ on L to the smooth (k+1)-form again on L $e^{-\tau\varphi} d^k_L(e^{\tau\varphi}\alpha)$.

We shall call Witten's tangential Laplacian the measurable filed $(\Delta^k_{\tau,L})_L$ which is defined in the obvious way. Then we prove that $\Delta^k_\tau$ computes the ($L^2$) tangential cohomology of (M,F):

*Proposition 1:* The fields (Ker $(\Delta^k_{\tau,L})_L$ and (Ker $(\Delta^k_L)_L$) of Hilbert spaces are measurably isomorphic and one has that

$$\beta_k = \dim_\Lambda (\text{Ker}(\Delta^k_{\tau,L})_L) < +\infty$$

for any positive real $\tau$ and $0 \leq k \leq p$.

*Proof (Sketch):* The proposition can be proved following the steps below:

1. First one has to prove that the operator $d^k_\tau = (d^k_{\tau,L})_L$ is a differential operator which is elliptic along the leaves of F. This can be proved using an argument similar to the one used by Connes in ([5]) to prove the transversal index theorem.

2. Next one proves that the adjoint of $d^k_{\tau,L}$ is the closure of the operator which sends each smooth (k+1)-form $\beta$ on L to the smooth k-form again on L $e^{\tau\varphi}(d^k_L)^*(e^{-\tau\varphi}\beta)$ where

$$(d^k_L)^* = (-1)^{pk+1} *d^k_L* ,$$

where * is the Hodge star operator on the leaf L defined via the Remanian metric.

3. We note that $\Delta^k_\tau = (\Delta^k_{\tau,L})_L$ is a field of measurable positive operators acting on the Hilbert space of square integrable k-forms on the leaf L. Moreover $\Delta^k_\tau$ is elliptic along the leaves.

4. For any leaf L and $0 \leq k \leq p$ we denote by $T^k_L$ the bounded operator on $L^2(L, \wedge^k T^*F)$ defined by

$$T^k_L(\alpha)(x) = e^{-\tau\varphi(x)} \alpha(x),$$

for $\alpha$ in $L^2(L, \wedge^k T^*F)$.

It is clear that $T^k_L$ is invertible and defines an element of $L^\infty(M/F, \wedge^k T^*F)$. Next we set

$$U^k{}_{,L} = Q^k{}_{,L}\, T^k{}_L\, Q^k{}_L$$

where $Q^k{}_{,L}$ (resp. $Q^k{}_L$) denotes the orthogonal projection onto the subspace Ker($\Delta^k{}_{,L}$) (resp. onto Ker($\Delta^k{}_L$)). We thus define a measurable field $(U^k{}_{,L})_L$ of endomorphisms of the random Hilbert space $(L^2(L, {}^kT^*F))_L$, such that Ker($\Delta^k{}_{,L}$) is a superset of $U^k{}_{,L}$(Ker($\Delta^k{}_L$)).

We want to show that $(U^k{}_{,L})_L$ belongs to L (M/F, ${}^kT^*F$) and defines an isomorphism of Hilbert spaces from (Ker($\Delta^k{}_L$))$_L$ to (Ker($\Delta^k{}_{,L}$))$_L$. One then has (omitting the subscript L):

$$d^k{}_\varepsilon = T^{k+1}{}_\varepsilon\, d^k (T^k{}_\varepsilon)^{-1}$$

and hence

$$T^k{}_\varepsilon (\text{Ker } d^k) = \text{Ker } d^k{}_\varepsilon \quad \textbf{(1)}$$

and

$$T^{k+1}{}_\varepsilon (\text{cl.Im}(d^k)) = \text{cl.Im}(d^k{}_\varepsilon). \quad \textbf{(2)}$$

But it follows from Hodge theory that one has the following orthogonal decompositions:

$$\text{Ker}(d^k) = \text{Ker}(\Delta^k) + \text{cl.Im}(d^{k-1})$$

and

$$\text{Ker}(d^k{}_\varepsilon) = \text{Ker}(\Delta^k{}_\varepsilon) + \text{cl.Im}(d^{k-1}{}_\varepsilon),$$

and then from equations (1) and (2) it follows that $T^k{}_\varepsilon$ is given in those decompositions by a 2 X 2 matrix with the upper left entry being $U^k$, the lower right entry being $B^k$, the upper right entry being 0 and the lower left entry being any element and where the entry $B^k{}_{,L} = \Delta^k{}_{,L} \mid \text{cl.Im}(d^{k-1}{}_L)$, namely $\Delta^k{}_{,L}$ restricted to cl.Im($d^{k-1}{}_L$), the closure of the Image of $d^{k-1}{}_L$ is invertible. We thus deduce that $U^k$ is an isomorphism from (Ker($\Delta^k{}_L$))$_L$ onto (Ker($\Delta^k{}_{,L}$))$_L$ and hence

$$\beta_k = \dim_\Lambda (\text{Ker}(\Delta^k{}_L)_L) = \dim_\Lambda (\text{Ker}(\Delta^k{}_{,L})_L)$$

and this holds for all $\varepsilon > 0$. As $\beta_k < +\infty$ then an argument similar to Connes [6] completes the proof.